\input amstex
\documentstyle{amsppt}
\magnification=1200
\pagewidth{6.5 true in}
\pageheight{9.0 true in}
\nologo
\tolerance=10000
\overfullrule=0pt
\define\m{\bold m}
\define\A{\frak A}
\define\a{\underline{a}}
\define\la{\longrightarrow}
\define\inc{\subseteq}
\define\kd{\text{dim}}
 
\topmatter
\title
A cancellation theorem for ideals
\endtitle
\rightheadtext{normal ideals}
\author
Craig Huneke
\endauthor
\address
Department of Mathematics, University of Kansas,
Lawrence, KS 66045
\endaddress
\email
huneke\@math.ukans.edu
\endemail
\abstract
We prove  cancellation theorems for special ideals in Gorenstein local
rings. These theorems take the form that if $KI\inc JI$, then $K\inc J$.
\endabstract
\thanks
The author is grateful to the Max-Planck Institue for
its support while visiting during Fall 1998,
where part of this work was completed. The author also
thanks the NSF for partial support. 
\endthanks
\endtopmatter
 
\document
 
\head 1. Introduction
\endhead
One of the basic themes in commutative Noetherian ring
theory is that of cancellation of ideals: if $I,J,K$ are
ideals in a commutative ring $R$ and $IK = JK$, then one
cannot in general conclude that $I = J$. Of course for
Dedekind domains it is an old and famous theorem that
one can do so. In general, one must
be satisfied with `cancellation up to integral closures',
i.e. in concluding that the integral closure of $I$ is
equal to the integral closure of $J$.

It is of considerable interest to understand when one can cancel: rarely 
is what experience teaches. However, the point of this paper is to
prove a cancellation theorem for special ideals. This cancellation has
a number of immediate corollaries, e.g. to theorems of Brian\c con-Skoda type
which say that the integral closure of a certain power of $I$  is contained
in $I$, as well as  a number of corollaries dealing with reduction
numbers and the syzygetic property.

We recall relevant definitions:

\definition{Definition 1.1} Let $R$ be a ring and let $I$ be an ideal of $R$.
 An element $x\in R$ is
{\it integral over $I$\/} if $x$ satisfies an equation of the form $x^n+a_1x^{
n-1}+\dots+a_n=0$, where $a_j\in I^j$ for $1\leq j\leq n$.
The \it integral closure \rm of $I$, denoted by $\overline I$,
is the set of all elements integral over $I$. This set is an ideal.
\enddefinition

Closely related to integral closure is the idea of a reduction.
Recall that $J\inc I$ is
said to be  a \it reduction \rm of $I$ if there exists an integer $n \geq 0$
such that $I^{n+1} = JI^{n}$. The least such $n$ is said to be the reduction
number of $I$ with respect to $J$.  It is denoted $r_J(I)$. If $J$ is a reduction
of $I$ and properly contains no other reduction of $I$, $J$ is said to be
a \it minimal reduction\rm.

The \it analytic spread \rm of an ideal $I$ in a local Noetherian ring $(R,\m)$
with infinite residue field is the least number of generators of any
minimal reduction of $I$. Equivalently, it is the Krull dimension of
the `fiber cone', $R[It]/\m R[It]$. The \it analytic deviation \rm of $I$,
denoted ad$(I)$,
is the analytic spread of $I$ minus the height of $I$.

Our main theorem is Theorem 2.2 below. Its most basic case says that
if $R$ is a regular local ring, $P$ a prime of $R$ with dim$(R/P) = 1$,
and $J$ any minimal reduction of $P$, then $IP\inc JP$ implies $I\inc J$
for an arbitrary ideal $I$. 
\bigskip
\head 2. Main Result\endhead

We begin with an elementary lemma.

\proclaim{Lemma 2.1} Let $(R,\m)$ be a local Gorenstein ring. Let
$I\inc R$ be an ideal  of height $g$. Set $d = \text{dim}(R/I)$.
Let $\a = (a_1,...,a_g)$ be a regular sequence in
$I$. 
Then $R/(\a:I)$ satisfies Serre's condition
$S_{2}$ iff $H^{d-1}_{\m}(R/I) = 0$.
\endproclaim

\demo{Proof} We reduce first to the case in which $I$ is unmixed.
Let $K$ be the unmixed part of the primary decomposition of $I$,
so that $K/I$ has dimension strictly less than $d$. The long
exact sequence on local cohomology applied to the short
exact sequence, $$0\la K/I\la R/I\la R/K\la 0$$ shows that
$H^{d-1}_{\m}(R/K) = 0$ as it lies between $H^{d-1}_{\m}(R/I)$, which
is $0$ by assumption, and $H^d_{\m}(K/I)$, which is $0$ since
$d$ is greater than the dimension of $K/I$. We can replace $I$ by
$K$ if we prove that $(\a:I) = (\a:K)$. It suffices to check this
at the associated primes of $(\a:K)$ (as $I\inc K$ implies that
$(\a:K)\inc (\a:I)$),
and all of these are assocated to $\a$. If $Q$ is such a prime and
 $I\nsubseteq Q$, then $(\a)_Q = (\a:I)_Q = (\a:K)_Q$. If $I\inc Q$,
then $I_Q = K_Q$, and again $(\a:I)_Q = (\a:K)_Q$. We replace $I$ by $K$
and assume it is unmixed.

Assume that $H^{d-1}_{\m}(R/I) = 0$.
Let $Q$ be a prime containing $\a:I$, and assume that the
height of $Q$ is at least $g+1$. 
If $Q$ does not contain $I$, then
$(\a:I)_Q = (\a)_Q$ and the result follows. Henceforth we assume $I\inc Q$.

Choose an element $t\in Q$ which is not a zero-divisor on either
$R/I$ or $R/(a_1,...,a_g)$. Apply Hom$_R(R/I,\quad)$ to the
short exact sequence $0\la R\overset {t}\to\la R\la S\la 0$,
where $S = R/Rt$. Note that Ext$^{g+1}_R(R/I,R)$ is the Matlis
dual of $H^{d-1}_{\m}(R/I)$ (recalling $R$ is Gorenstein), and
hence is zero. Thus we get a short exact sequence,
$$0\la \text{Ext}^g_R(R/I,R)\overset{t}\to\la \text{Ext}^g_R(R/I,R)\la \text{Ext}^g_R(R/I,S)\la 0.$$
Since $a_1,..,a_g$ are a regular sequence on $R$ and $R/tR$ in $I$,
we may identify
this last exact sequence with the sequence,
$$0\la \text{Hom}_R(R/I,R/\underline{a})\overset{t}\to\la \text{Hom}_R(R/I,R/\underline{a})\la \text{Hom}_R(R/I,S/\underline{a}S)\la 0.$$
This gives the identity,
$$ (\a:I) + tR = (\a+Rt:I). \tag{1}$$
In particular, since $a_1,...,a_g,t$ form a regular sequence, the ideal
$(\a:I) + tR$ is unmixed, necessarily of height $g+1$. But then
depth$((R/((a_1,...,a_g):I))_Q) \geq 2$, unless $Q$ has height exactly $g+1$,
in which case this localization is Cohen-Macaulay. Thus $R/(\a:I)$ satisfies
$S_2$.

Conversely, suppose that $R/(\a:I)$ satisfies 
$S_2$. 
Notice that Ext$^{g+1}_R(R/I,R)$ has an annihilator of height at least $g+2$
since at all primes of height $g+1$, after localization this module is
zero as $R/I$ is $S_1$. Choose an element $t$, a zero-divisor on both
$R/I$ or $R/(a_1,...,a_g)$ such that $t$Ext$^{g+1}_R(R/I,R) = 0$. The same
exact sequences as above then give the isomorphism
$$ ((\a+Rt):I)/((\a:I) + tR)\cong \text{Ext}^{g+1}_R(R/I,R).$$
This proves that the associated primes of $ \text{Ext}^{g+1}_R(R/I,R)$ are
contained in the associated primes of $R/((\a:I) + tR)$, and since
$R/(\a:I)$ satisfies  
$S_2$, all such primes have height at most $g+1$. It follows that
$ \text{Ext}^{g+1}_R(R/I,R) = 0$. \qed\enddemo

\proclaim{Theorem 2.2 (Cancellation)} Let $(R,\m)$ be a local Gorenstein ring. Let
$I\inc R$ be an unmixed ideal  of height $g$,
and let $a_1,...,a_{g+1}\in I$. Assume that $a_1,...,a_g$ form a regular
sequence which  generates
$I$ locally at the localization of every minimal of $I$, and that
$a_{g+1}$ is not contained in any minimal prime of $(a_1,...,a_g)$
which does not contain $I$. Set $J = (a_1,...,a_g,a_{g+1})$.
Finally assume that $H^{d-1}_{\m}(R/I) = 0$ (e.g. if $R/I$ is Cohen-Macaulay),
 where $d = \kd(R/I)$.
If $KI\inc JI$, then $K\inc J$.
\endproclaim

\demo{Proof}
Observe that the assumptions of the theorem force $\a:I = \a:a_{g+1}$.
We first reduce to the situation in which $R/I$ has dimension one;
this is the main case. Our assumption on the local cohomology is
made to reduce to this case. 

To prove $K\inc J$, it suffices to prove $K_Q\inc J_Q$ for every associated
prime $Q$ of $J$. We claim all such associated primes have height at
most $g + 1$. This will then reduce to the case in which $R/I$ has
dimension one by localization. Let $Q$ be an associated prime
of $J$ and suppose that $Q$ has height at least $g + 2$. 
The exact sequence
$$0\la R/((a_1,...,a_g):a_{g+1})\la R/(a_1,...,a_g)\la R/J\la 0$$
induced by multiplication by $a_{g+1}$ on the first map shows that
$Q$ must be associated either to $R/(a_1,...,a_g)$, which is
impossible, or depth$((R/((a_1,...,a_g):a_{g+1}))_Q) \leq 1$.
By Lemma 2.1, $R/((a_1,...,a_g):a_{g+1})$ is $S_2$. As $Q$ has height at least $g + 2$,
this is a contradiction, proving our claim.

We have reduced to the case in which dim$(R/I) = 1$. Choose an element
$s$, not a zero-divisor on $R/I$, such that $\a:s = I$. This is possible
since $I$ is generically equal to $\a$ and is unmixed: specifically choose
an element $s$ which is in all the primary components of the ideal $\a$ whose
radicals do not contain $I$ and in addition choose $s$ not in any associated
prime of $I$. As both $I$ and $\a$ are unmixed this is possible using
prime avoidance. The condition that $s$ not be a zero-divisor on $R/I$ is
immediate from the choice made. We can write $\a = N\cap I$, where the minimal
primes over $N$ do not contain $I$. 
We have chosen $s\in N$ but not a zero-divisor modulo $I$. Hence
$\a:s = (N:s) \cap (I:s) = I$.
Set $b = a_{g+1} + s$,
and set $\frak A = (a_1,...,a_g,b)$.
Observe that $\A$ is $\m$-primary, as $b$ is not contained in the union
of the minimal primes over $\a$.

Set $q = (\a:I) + I$. We claim that $K\inc \A:q$. Clearly $K(\a:I)\inc \a\inc \A$,
since $K\inc I$. This last containment follows since it is enough to check
it at the localizations of the associated primes of $I$, which are all minimal
over $I$. After localizing at such a prime, $I$ is equal fo $\a$, and our
assumption then reads (after localization at a minimal prime $Q$) that
$K_Q(\a)_Q\inc (\a)_Q^2$. As $\a$ is generated by a regular sequence it then
follows that $K_Q\inc (\a)_Q$. Hence $K\inc I$.  To prove that $KI\inc \A$,
it is enough to prove that $JI\inc \A$ by our basic assumption. As $\a I\inc \A$,
we only need to prove that $a_{g+1}I\inc \A$. But for $i\in I$, we may write
$a_{g+1}i = (a_{g+1}+s)i - si\in \A$ since $si\in \a$.

We next claim that $\A:s = q$. One direction is much the same as above.
We know $sI\inc \a$ by choice of $s$. If $t\in \a:I$, then we can write
$ts = tb - ta_{g+1}\in \A$ since $ta_{g+1}\in \a$. Thus $q\inc \A:s$.
Conversely, suppose that $su\in \A$. Write $su = bv + w$, where $w\in \a$.
Then $s(u-v)\in J\inc I$, and thus $u-v\in I$.  Note that $Ibv\inc Isu + Iw\inc
\a$. Since $b$ is a regular element modulo $\a$, we see that $Iv\inc \a$.
Thus $u = (u-v) + v\in I + \a:I = q$.

As $\A$ is an $\m$-primary ideal in a Gorenstein local ring, generated by
a regular sequence, we know that $K\inc \A:q = \A:\A:s = \A + Rs = J + Rs$.
Hence $K\inc (J + Rs)\cap I = J + (I\cap Rs)\inc J + Is\inc J$, finishing 
the proof. \qed\enddemo

\smallskip

\proclaim{Corollary 2.3}  Let $(R,\m)$ be a local Gorenstein ring with
infinite residue field and
suppose that $I$ is an unmixed ideal, generically a complete intersection,
such that the ad$(I) = 1$. Assume that $H^{d-1}_{\m}(R/I) = 0$ where
$d$ is the Krull dimension of $R/I$. Let $J$ be an arbitrary minimal
reduction of $I$. Then $r_J(I)\ne 1$.
In particular, if $R$ is a regular ring and $P$ is a prime with dim$(R/P) = 1$,
then $P$ cannot have reduction number $1$ with respect to a minimal
reduction of $P$.
\endproclaim

\demo{Proof} The last statement follows immediately from the first.
The assumptions guarantee  that a given minimal reduction $J$ of $I$ is
generated by $g+1$ elements which can chosen as in the conditions of Theorem
2.2. If the reduction number is one, then $P^2 = JP$, and cancellation
then proves that $P\inc J$, that is $P = J$, and the reduction number is $0$.\qed
\enddemo

\proclaim{Theorem 2.4} Let $(R,\m)$ be a regular local ring with infinite
residue field. Suppose that $I$ is an unmixed ideal,
generically a complete intersection, 
such that the ad$(I) = 1$. Assume that $H^{d-1}_{\m}(R/I) = 0$ where
$d$ is the Krull dimension of $R/I$. Set $g = \text{height}(I)$. Then
$\overline{I^g}\inc J$, where $J$ is an arbitrary minimal reduction of $I$.
In particular, if $P$ is a prime of $R$ with dim$(R/P) = 1$, then
$\overline{P^g}\inc J$, for an arbitrary minimal reduction $J$ of $P$. 
\endproclaim

\demo{Proof} By the results in \cite{AH} and \cite{L}, we know that
$$\overline{I^{g+1}}\inc JI.$$ 
Since $I(\overline{I^g})\inc \overline{I^{g+1}}$, the result follows
at once from the cancellation theorem. \qed\enddemo

The above result suggests at least the possibility that for every prime $P$
in a regular local ring $R$, $\overline{P^g}\inc J$, for $g = $ height$(P)$
and $J$ an arbitrary minimal reduction of $P$. This would be a considerable
improvement on the usual Brian\c con-Skoda results which only give
$\overline{P^d}\inc J$ in general, for $d = \kd (R)$. Trying to find an
example in which the stronger claim $\overline{P^g}\inc J$ does not
hold turned out to be difficult, and no counterexample was found. In fact,
nearly every example checked, no matter what the analytic spread of $P$,
had the containment $P^2\inc J$. Example 2.7 below gives a height three
prime ideal such that $P^2\nsubseteq J$. The examples below deal with
powers, rather than the integral closure of powers as this author does
not know a good way to compute the integral closures of the powers of
these ideals. The examples were all computed on MACAULAY, in the highest
characteristic allowable.

\example{Example 2.5} Let $P$ be the defining ideal of
$k[t^4+s^4,t^3s^2,t^2s^4+ts^5,ts^6]$ in the ring $R = k[X_1,...,X_4]$.
Then $P$ is height two,  is generated by $5$ elements, say $a_1,...,a_5$,
and $R/P$ is not Cohen-Macaulay. The fiber cone
$S/\m S$, where $\m = (X_1,...,X_4)$ and $S = R[Pt]$, is isomorphic to
$k[T_1,...,T_5]/(T_1T_5)$. It follows that the analytic spread of 
$P$ is $4$. A minimal reduction is $J = (a_1+a_5,a_2,a_3,a_4)$.
One can check that $P^2\inc J$.
\endexample

\example{Example 2.6} For a slightly more complicated example, consider the
defining ideal of the ring $k[s^3,t^3,u^3,s^2t+stu,st^2+tu^2]$. This ideal
$P$ is height two, generated by $8$ elements.  If we map $R[A,B,C,D,E,F,G,H]$
onto the Rees algebra using these $8$ generators, the fiber cone
$k[A,B,C,D,E,F,G,H]/I$ is defined by $I = (F^2-1/2DG+EG-FG, EF-1/2DG-2EG-FG+G^2,
E^2+4/9CG+5/6DG+EG+FG-G2+1/18EH-1/18FH-1/6GH DF-CG,
DE-1/9CG-1/3DG+1/9EH-1/9FH-1/3GH, D^2+16/9CG-2/3DG+2/9EH-2/9FH+4/3GH, CF+GH,
CE-1/3CG+1/3EH-1/3FH, CD-2/3CG+2/3EH+4/3FH-2GH, C^2+DH)$. This ideal has height
four, so the analytic spread of $P$ is $4 = 8-4$. A minimal reduction $J$ is
given by the first three generators corresponding to $A,B,C$, together with
the sum of all the rest of the specified generators. Precisely, if we map
$k[a,b,c,d,e]$ onto $k[s^3,t^3,u^3,s^2t+stu,st^2+tu^2]$, then  $J$ is
generated by the four elements,

\noindent $f_1 = a^2bc-bc^3+3abcd-3bc^2d+3abd^2-3bcd^2-bd^3-cd^3-3a^2be-3abce$
\newline\noindent $-3abde-3d^3e+3ade^2 +3d^2e^2+ae^3+ce^3$,
\smallskip
\noindent $f_2 = a^2b^2+1/2abc^2-1/2bc^3-3/2bc^2d-3/2bcd^2-1/2bd^3-3abce-3/2abde$
\newline\noindent$ +3/2d^2e^2- 1/2ae^3+1/2ce^3$,
\smallskip
\noindent $f_3 = abd^3+c^2d^3-a^2bde-abcde-3abd^2e-2d^4e-a^2be^2-abce^2-3acde^2$
\newline\noindent$ +3ad^2e^2+d^3e^2 -a^2e^3-ace^3+3ade^3$, and
\smallskip
\noindent $f_4 = ab^2c^2+b^2c^3+ab^2cd+abc^2d+3b^2c62d+2bc^3d+a^2bd^2+2ab^2d^2+abcd^2$
\newline\noindent$ +3b^2cd^2 +5bc^2d^2+b^2d^3+3bcd^3-3c^2d^3+bd^4-3cd^4-d^5-a^3be-3ab^2ce+4abc^2e$
\newline\noindent$ -4a^2bde -7ab^2de+2abcde+2bc^2de+ad^3e-3bd^3e-cd^3e+d^4e+3a^2be^2+2ab^2e^2$
\newline\noindent$ -6abce^2+2bc^2e^2 +3abde^2+9acde^2+6bcde^2+3bd^2e^2+3a^2e^3+5abe^3+3ace^3$
\newline\noindent$ -bce^3-4ade^3-2cde^3+d^2e^3 -3ae^4-2de^4-2e^5.$
 
One can check that $P^2\inc J$.
\endexample

\example{Example 2.7} Let $I$ be the ideal of maximal minors of the 
$4$ by $6$ matrix,
$$
A = \left( \matrix
a & b & c & e & f & d\\
d & a & b & c & e & f\\
e & d & a & b & c & e\\
f & e & d & a & b & c
\endmatrix \right).
$$

The ideal $I$ has generic height in $R = k[a,b,c,d,e,f]$, namely $3 = 6-4+1$. 
In particular $R/I$ is Cohen-Macaulay.
This ideal was chosen so that after the specialization $ d = e = f = 0$,
the specialized ideal becomes the $4$th power of $(a,b,c)$, and the square
of this ideal is not contained in any minimal reduction of it.  One can check
(with a huge amount of computer time!) that $I$ is reduced, and has analytic
spread $6$.   We computed a particular minimal reduction $J$ (which is
extremely complicated), and for this particular reduction, $P^2\nsubseteq J$.
However, $P^2\inc J:\m$, where $\m$ is the obvious maximal ideal. In
particular $P^3\inc J$. The calculation of this example took more than a
week of computer time, and discouraged the author from looking much further!
\endexample

We recall some concepts before stating the next theorem. Let $R$ be a
Noetherian ring and $I$ an ideal. If $I$ is generated by $a_1,...,a_n$
we can map $S = R[T_1,...,T_n]$ onto $R[It]$ by sending $T_i$ to $a_it$. The
kernel of this map we denote by $Q$. This ideal is graded, and by
$Q_j$ we denote the $j$th graded piece of $Q$. The ideal $I$ is said
to by \it syzygetic \rm if $Q_2\inc S_1Q_1$. See \cite{V, 2.1} for a discussion
of this concept. Some known classes of syzygetic ideals are of course
those of linear type, where $SQ_1 = Q$, but also include, e.g., all
height two reduced ideals in a $3$-dimensional regular local ring.

\proclaim{Theorem 2.8} Let $R$ be a Gorenstein local ring and $P$ a prime
ideal such that ad$(P) = 1$, $R_P$ is regular, and $H^{d-1}_{\m}(R/P) = 0$ where
$d$ is the Krull dimension of $R/P$. If $\mu (P) = g+2$, where $g =$ height$(P)$,
then $P$ is syzygetic.
\endproclaim

\demo{Proof} Suppose not. Fix a generating set $a_1,...,a_{g+2}$ of $P$.
 Then there is a homogeneous polynomial of
degree two, $H\in R[T_1,...,T_{g+2}]$ such that $H(a_1,...,a_{g+2}) = 0$
and such that $H\notin \frak A$, where $\frak A = SQ_1$ is the
ideal defining the symmetric
algebra of $P$, i.e. the ideal in $R[T_1,...,T_{g+2}]$ generated by
all linear polynomials which vanish at $a_1,...,a_{g+2}$. After possibly
changing the generators, we may assume $H = rT_1^2 + P$, where
$P$ is a polynomial in the ideal $(T_2,...,T_{g+2})$; we may further
assume that $J = (a_2,...,a_{g+2})$ is a minimal reduction of $P$, and
that $P$ is generically generated by $a_2,...,a_{g+1}$ which form a
regular sequence. But
then $ra_1^2\in JI$, and hence $rI^2\inc JI$. By the cancellation theorem,
$rI\inc J$, so there is a linear polynomial $L = rT_1 + L'$, where
$L'\in (T_2,...,T_{g+2})$, such that $L(a_1,...,a_{g+2}) = 0$.
Consider the polynomial $H - T_1L = G$. The polynomial $G$ can be written
in the form $$G = T_2L_2 + ... + T_{g+2}L_{g+2}$$
for linear polynomials $L_i$. Set $l_i = L_i(a_1,...,a_{g+2})$. Then
$l_{g+2}\in (a_2,...,a_{g+1}):a_{g+2}\cap P$. Since $a_2,...,a_{g+1}$
generically generate $P$ this intersection is just $(a_2,...,a_{g+1})$.
Hence there is another linear polynomial $L* = r_2T_2+...+r_{g+1}L_{g+1}$
such that $l_{g+2} = L*(a_2,...,a_{g+1})$. In other words 
$L_{g+2} - L*$ is a linear polynomial in $\frak A$. Modulo $\frak A$  
we can write $G = T_2F_2 + ... + T_{g+1}F_{g+1}$ for some new choice of
linear polynomials $F_i$.  But then $G\in Q_2 \cap (T_2,...,T_{g+1})$,
and $a_2,...,a_{g+1}$ form a regular sequence. The next lemma then
finishes the proof.
It is well-known to experts, but a suitable reference
was not found. (See however \cite{AH, (2.5)}.)

\proclaim{Lemma 2.9} Let $(R,\m)$ be a Noetherian local ring and let
$I$ be an ideal. Assume that $ (a_1,...,a_n, a_{n+1},...,a_m) = I$, 
and that $\a = a_1,...,a_n$ form a regular sequence.
Let $R[T_1,...,T_m]\la R[It]$ be the
surjection of the polynomial ring $R[T_1,...,T_m]$ onto $R[It]$ sending
$T_i$ to $a_it$. Let $Q$ be the kernel of this surjection with the
obvious grading.
Then $(T_1,...,T_n)\cap Q_2 \inc S_1Q_1$.
\endproclaim

\demo{Proof} By descending induction on i, we prove that
$$(T_1,...,T_i)\cap Q_2 \inc S_1Q_1 + (T_1,...,T_{i-1})\cap Q_2.$$ 
Since $(T_1)\cap Q_2 = T_1Q_1$ as $a_1$ is a nonzerodivisor, the Lemma
will follow. 

Suppose that $G\in (T_1,...,T_i)\cap Q_2$, and write $G = \sum_{j=1}^i T_jG_j$,
with $G_j$ linear forms. Write $g_j = G_j(a_1,...,a_m)$. Then
$g_i\in (a_1,...,a_{i-1}):a_i = (a_1,...,a_{i-1})$, so we may write
$g_i = \sum_{k=1}^{i-1} a_kr_k$. Set $H = \sum_{k=1}^{i-1} T_kr_k$, so that
$G_i-H\in Q_1$. Then $G = \sum_{j=1}^{i-1} T_jG_j + T_i(G_i-H) + T_iH\in
(T_1,...,T_{i-1}) + S_1Q_1$, from which the result follows. \qed\enddemo
\enddemo

\proclaim{Theorem 2.10} Let $I$ be an unmixed ideal in a regular local ring of
dimension $d$.
such that $I_P$ is generated by a regular sequence of length $g < d$ for every
prime $I\inc P$ with $P\ne \m$. Then 
$\overline{I^{d-1}}\inc J$, for every reduction $J$ of $I$.
\endproclaim

\demo{Proof} We induce upon the dimension of $R/I$. If the dimension is $1$,
then Theorem 2.4 gives us the conclusion. Choose a regular sequence $\a = 
a_1,...,a_g\in I$ which generated $I_P$ for every minimal prime $P$ above $I$.
Set $K = \a:I + I$. Then the dimension of $R/K$ is strictly less than the
dimension of $R/I$. We claim that $K$ is unmixed and 
$K_Q$ is generated by a regular sequence of length $g+1$ for every prime
$K\inc Q$ with $Q\ne \m$. This follows from the next Lemma:

\proclaim{Lemma 2.11} Let $I$ be an unmixed ideal in a local Cohen-Macaulay
ring $(R,\m)$, and assume that $I$ is generically generated by a
regular sequence of length $g$.
Let $\a = a_1,...,a_g$ be a regular sequence in
$I$ which generates $I$ locally at all minimal primes of $I$.
Define an ideal $K = (\a:I) + I$. Then $K$ is an unmixed ideal of height $g+1$.
If $I_P$ is generated by a
regular sequence of length $g$ for all primes $I\inc P$ such that
dim$((R/I)_P)\leq 1$, then $K$ is generically a complete intersection.
More generally if $K\inc Q$
and $I_Q$ is generated by a regular sequence, then $K_Q$ is also
generated by a regular sequence.
\endproclaim
 
\demo{Proof} We first prove that the height of $K$ is at least $g+1$.
If $P$ contains $K$ and has height at most $g$, then $P$ must contain $I$
(hence $P$ has height exactly $g$) and $(\a)_P = I_P$ by assumption.
Then $\a:I\nsubseteq P$, a contradiction. Hence height$(K)\geq g+1$.
Let $Q$ be an associated prime of $K$. The exact sequence
$$0\la R/((\a:I)\cap I)\la R/I\oplus R/(\a:I)\la R/K\la 0$$
shows that either $Q$ is associated to $I$, to $\a:I$, or else
depth$(R/((\a:I)\cap I))\leq 1$. Both $I$ and $\a$ are unmixed, so
every associated prime of either $I$ or $\a:I$ has height $g$. On the other
hand, since $K((\a:I)\cap I)\inc \a$, and height$(K)\geq 1$, it follows
that $(\a:I)\cap I = \a$. Then the height of $Q$ is at most $g+1$, which
proves that $K$ is unmixed of height $g+1$.
 
The penultimate statement follows at once from the more general last statement
of the theorem.
Let $Q$ be a prime containing $K$, such that  $I_Q$ is a
complete intersection, say generated by $b_1,...,b_g$. Then $\a:(b_1,...,b_g)$
is generated by $a_1,...,a_g, \delta$, where $\delta$ is the determinant of
a transition matrix $D$ writing the $a_i$ in terms of the $b_i$. See
\cite{BH,(2.3.10)}. It follows that $K_Q = (a_1,...,a_g, \delta)_Q + (b_1,...,b_g)_Q =
(b_1,...,b_g, \delta)_Q$ is generated by $g+1$ elements, which are necessarily
a regular sequence as the height of $K$ is $g+1$ and $R$ is Cohen-Macaulay.
\enddemo

We can apply the induction to conclude that $\overline{K^{d-1}}\inc L$,
for every reduction $L$ of $K$. Let $J$ be an arbitrary reduction of $I$.
Then $J + \a:I$ is certainly a reduction of $K$, so we obtain that
$$\overline{I^{d-1}}\inc \overline{K^{d-1}}\cap I\inc (J+\a:I)\cap I = J + (\a:I\cap I) = J$$
which finishes the proof.\qed\enddemo

Another case is easy:

\proclaim{Theorem 2.12} Let $I$ be an ideal in a regular local ring $(R,\m)$
 of dimension $d$. Assume that $I$ is generically a complete intersection, and
that $R/I$ is Cohen-Macaulay of positive dimension. Then $\overline{I^{d-1}}\inc J$ for every reduction
$J$ of $I$.
\endproclaim

\demo{Proof} Induce upon the dimension of $R/I$. If this dimension is $1$, we
obtain the conclusion from Theorem 2.4. Assume it is bigger than $1$. Denote
the height of $I$ by $g$.The set of
primes $q$ of height $g+1$ such that $I_q$ is not a complete intersection is finite
in number, say $q_1,...,q_n$. Choose 
an element $y\in \m$ a nonzerodivisor on $R/I$, and such that $y$ is not
contained in the union of $q_1,...,q_n$. Such a choice is possible by prime
avoidance. Set $K = I + Ry$. Then
$K$ is generically a complete intersection, since if $Q$ is minimal over $K$
then $Q\ne q_i$ for any $i$, and hence $I_Q$ is generated by $g$ elements and
$K_Q$ by $g+1$ elements. Clearly $R/K$ is Cohen-Macaulay. Fix a reduction $J$
of $I$. Then $J+Ry$ is a reduction of $K$.  The induction 
implies that $$\overline{I^{d-1}}\inc \overline{K^{d-1}}\cap I\inc (J+Ry).$$
We may replace $y$ by $y^n$ and use the Krull Intersection Theorem to prove
that $\overline{I^{d-1}}\inc J$.\qed\enddemo
\bigskip

We end by giving one other somewhat surprising corollary of our main
theorem.

\proclaim{Corollary 2.13}
Let $(R,\m)$ be a local Gorenstein ring. Let
$I\inc R$ be an unmixed ideal  of height $g$,
and let $a_1,...,a_{g+1}\in I$. Assume that $a_1,...,a_g$ form a regular
sequence which  generates
$I$ locally at the localization of every minimal of $I$, and that
$a_{g+1}$ is not contained in any minimal prime of $(a_1,...,a_g)$
which does not contain $I$. Set $J = (a_1,...,a_g,a_{g+1})$.
Finally assume that $H^{d-1}_{\m}(R/I) = 0$ (e.g. if $R/I$ is Cohen-Macaulay),
 where $d = \kd(R/I)$ and that the radical of $J$ is equal to
the radical of $I$ (e.g. if $J$ is a reduction of $I$). Then
$I^n\inc J$ iff $I^{n+1}\cap J\inc JI$.
\endproclaim

\demo{Proof} One direction is obvious. If $I^n\inc J$, then $I^{n+1}\inc
I^n\cdot I\inc JI$. For the other direction we assume that $I^{n+1}\cap J\inc JI$
and prove by descending induction on $N$ that $I^N\inc J$. When $N = n$ we are
done. By assumption, some power of $I$ is contained in $J$. Suppose that
$I^N\inc J$ and $N > n$. Then $I^N\inc I^{n+1}\cap J\inc JI$, and by Theorem
2.2, we may cancel the $I$ to obtain that $I^{N-1}\inc J$. \qed\enddemo

\remark{Acknowledgement} The author thanks Ian Aberbach for many valuable conversations
concerning the material in this paper.
\endremark

\bigskip
\centerline{\bf Bibliography}
\bigskip
\refstyle{A}


\Refs\nofrills{}
\widestnumber\key{HH12}

\ref
\key {AH}
\by I. M. Aberbach and C. Huneke
\paper An improved
Brian\c con-Skoda theorem with applications to the Cohen-Macaulayness of Rees
rings
\jour Math. Ann.
\vol 297 \yr 1993
\pages 343--369
\endref

\ref
\key{BH}
\by W. Bruns and J. Herzog
\book Cohen-Macaulay Rings
\publ Cambridge studies in advanced mathematics
\vol 39
\yr 1993
\endref

\ref
\key L
\by J. Lipman
\paper Adjoints of ideals in regular local rings
\jour Math. Res. Letters
\vol 1
\yr 1994
\pages 1--17
\endref

\ref\key{V}
\by W. V. Vasconcelos
\book Arithmetic of Blowup Algebras
\publ Cambridge University Press
\publaddr Cambridge \yr 1994
\endref

\endRefs
\enddocument